\documentclass[authoryear]{elsarticle}

\usepackage{amsmath}
\usepackage{amssymb}
\usepackage{amsfonts}
\usepackage{mathptmx} 
\usepackage{bm}
\usepackage{tikz}
\usepackage{pgfplots}
\usepackage{color, colortbl}
\usepackage{graphicx}

\graphicspath{{pictures/}}

\usepackage{hyperref}

\newtheorem{thm}{Theorem}

\newproof{pf}{Proof}

\journal{arXiv.org}

\begin{document}

\begin{frontmatter}

\title{Two-component domain decomposition scheme with overlapping subdomains for parabolic equations}

\author[nsi,univ]{Petr N. Vabishchevich\corref{cor}}
\ead{vabishchevich@gmail.com}

\address[nsi]{Nuclear Safety Institute, Russian Academy of Sciences, 52, B. Tulskaya, Moscow, Russia}
\address[univ]{North-Eastern Federal University, 58, Belinskogo, Yakutsk, Russia}

\cortext[cor]{Corresponding author}

\begin{abstract}
An iteration-free method of domain decomposition is considered for approximate solving a boundary value problem for a second-order parabolic equation. A standard approach to constructing domain decomposition schemes is based on a partition of unity for the domain under the consideration.
Here a new general approach is proposed for constructing domain decomposition schemes  with overlapping  subdomains based on indicator functions of subdomains.
The basic peculiarity of this method  is connected with a representation of the problem operator as the sum of two operators, which are constructed for two separate subdomains with the subtraction of the operator that is associated with the intersection of the subdomains. There is developed a two-component factorized scheme, which
can be treated as a generalization of the standard Alternating Direction Implicit (ADI) schemes
to the case of a special three-component splitting. There are obtained conditions for the unconditional stability of regionally additive schemes constructed using indicator functions of subdomains.
Numerical results are presented for a model two-dimensional problem.

\end{abstract}


\end{frontmatter}

\section{Introduction} 

Schemes of domain decomposition are considered for numerical solving time-dependent
problems to partial differential equations. Iteration-free algorithms of domain decomposition
take into account specific features of time-dependent problems in the most efficient way.
In some cases, it is possible (see, e.g., \cite{Kuznetsov1988,Kuznetsov1991})
without loss of accuracy of the approximate solution make only one
iteration of the Schwarz alternating method at a new time level in
solving boundary value problems for a parabolic second-order equation. 
Iteration-free schemes of domain decomposition are associated with various variants of additive schemes (splitting schemes) -- see regionally additive schemes in \cite{SamarskiiMatusVabischevich2002}.

Methods of domain decomposition for solving unsteady problems can be classified by
(I) the method of decomposition for a calculation domain, 
(ii) the choice of decomposition operators (exchange boundary conditions), and
(III) the splitting scheme (approximation in time) employed.
For multidimensional boundary value problems, it is possible to use domain decomposition methods
with or without overlapping subdomains \citep{QuarteroniValli1999,ToselliWidlund2005}.
Methods without overlapping  subdomains are connected with an explicit formulation of exchange data conditions
on subdomain interfaces.

To construct decomposition operators for solving unsteady initial boundary value (IBV) problems for partial differential equations, it is convenient to use a partition of unity 
for a computational domain \citep{Dryja1991,Laevsky1987,SamarskiiVabishchevich1995a,Vabischevich1989,Vabischevich1994b,Vabischevich1994a}.
In decomposition methods with overlapping  subdomains,  the functions are associated with individual subdomains and take a value between zero and one.
Results of studies on  domain decomposition methods for Cauchy problems for partial differential equations
are summarized in the book \cite{SamarskiiMatusVabischevich2002}.  Among more recent studies, we highlight the works  \cite{Vabishchevich2008,vab-2011},
where schemes of domain decomposition  more suitable for numerical implementation are presented.

The construction of regionally additive schemes and study of their convergence are carried out
on the basis of the general theory of splitting schemes \citep{Samarskii1989,Marchuk1990,Vabishchevich2014}.
We highlight the simplest case of two-component splitting. In this case, we obtain
unconditionally stable factorized splitting schemes, such as
classical methods of alternating directions, predictor-corrector schemes and so on.
A more interesting for computational practice is the situation, when a problem operator is divided into a sum of three or more noncommutative non-self-adjoint operators.
In the case of such a multicomponent representation, splitting schemes are constructed on the basis of the concept of summarized approximation. 
For parallel computers, additively averaged splitting schemes
are of particular interest. In the class of splitting schemes with full approximation \citep{Vabishchevich2014}, we highlight vector additive schemes,
when the original equation is transformed into a system of similar equations \citep{Abrashin1990,Vabischevich1996b,AbrashinVabischevich1998}. 
The most convenient approach for constructing additive operator--difference schemes
of multicomponent splitting is based on  regularization of difference schemes
\citep{Samarskii1967}, where stability is achieved via perturbations of operators of the difference scheme.

The regionally additive schemes constructed on the basis of a partition of unity have
certain defects. The most important among them is connected with the use of the uniquely
defined operators with the generating coefficients in the overlapping domains. This leads to
the fact that, for example, in the problems with  constant coefficients of the equations one needs to use an algorithm with
varying coefficients. In the paper \cite{Vabishchevich2008}, we constructed unconditionally
stable regionally additive schemes with overlapping subdomains, which are more comfortable for the
practical use than traditional schemes designed on the basis of a partition of unity.
The regularization principle for operator--difference schemes \citep{Samarskii1967}  makes possible to construct
schemes of component-wise splitting. In the present work, domain decomposition methods  with overlapping subdomains
are designed using a two-component splitting, which generalize standard factorized schemes  \citep{SamarskiiMatusVabischevich2002,Vabishchevich2014}.
Decomposition operators are constructed on the basis of  the indicator functions of the subdomains.

The paper is organized as follows. A model problem for a parabolic equation with self-adjoint
elliptic operator of second order is formulated in Section 2. Approximation in space is constructed using Lagrangian finite elements,
whereas approximation in time is based on conventional two-level schemes with weights.
In Section 3, a two-component scheme of domain decomposition with overlapping  subdomains is constructed 
on the basis of a partition of unity for the computational domain.
New schemes of domain decomposition developed using indicator functions of subdomains are proposed in Section 4.
They are based on the generalization of classical factorized schemes (operator analogs of ADI schemes).
In Section 5, numerical experiments on the accuracy of the domain decomposition schemes are discussed for
the model IBV problem. The results of the work are summarized in Section 6.

\section{Problem formulation}

Let $\Omega$ be a bounded domain ($\Omega \subset \mathbb{R}^m, m = 2,3$) with a piecewise smooth boundary
$\partial \Omega$. We define an elliptic operator $\mathcal{A}$ so that
\begin{equation}\label{1}
 \mathcal{A} u = - \nabla (k(\bm x) \nabla u ) + c(\bm x) u ,
 \quad \bm x \in \Omega  
\end{equation} 
on the set of functions
\begin{equation}\label{2}
 u (\bm x) = 0 , 
 \quad \bm x \in \partial \Omega . 
\end{equation} 
Let  $k(\bm x)$ and $c(\bm x)$ be smooth functions in $\overline{\Omega}$ and
\[
 k(\bm x) \geq \kappa > 0,
 \quad c(\bm x) \geq 0 ,
 \quad \bm x \in \Omega . 
\] 
The Cauchy  problem
\begin{equation}\label{3}
 \frac{d u}{d t} + \mathcal{A} u = f(t),
 \quad 0 < t \leq T ,  
\end{equation} 
\begin{equation}\label{4}
 u(0) = u^0,
 \end{equation} 
is considered with $f(\bm x,t) \in L_2(\Omega)$, $u^0(\bm x) \in L_2(\Omega)$ using 
the notation $u(t) = u(\bm x,t)$. 

Let $(\cdot,\cdot), \|\cdot\|$ be the scalar product and norm in $H = L_2(\Omega)$, respectively:
\[
(u,v) = \int_\Omega u({\bm x})v({\bm x}) d {\bm x},
\quad \|u\| = (u,u)^{1/2}.
\]
A symmetric positive definite bilinear form $d(u,v)$ such that
\[
d(u,v) = d(v,u), 
\quad  d(u,u) \geq \delta \|u\|^2,
\quad  \delta > 0,
\]
is associated with the Hilbert space $H_d$, where the scalar product and norm are, respectively:
\[
(u,v)_d  = d(u,v), \quad  \|u\|_d  = (d(u,u))^{1/2}.
\]

Define $H^1_0(\Omega)$ as a subspace of $H^1(\Omega)$ such that
\[
 H^1_0(\Omega) = \{ v(\bm x) \in H^1(\Omega): \ v(\bm x) = 0, 
 \quad \bm x \in \partial \Omega \}.
\]  
Multiplying (\ref{3}) by $v(\bm x) \in H^1_0(\Omega)$ and integrating over the domain $\Omega$,
we arrive at the equality
\begin{equation}\label{5}
 \left( \frac{d u}{d t}, v \right ) + a(u,v) = (f, v), 
 \quad \forall  v \in H^1_0(\Omega) ,
 \quad 0 < t \leq T .
\end{equation} 
Here $a(\cdot, \cdot)$ is the following bilinear form:
\[
 a(u,v) = \int_{\Omega} ( k \nabla u \cdot \nabla v + c \, u \, v)  d \bm x ,
\] 
where
\[
 a(u,v) = a(v,u),
 \quad a(u,u) \geq \delta \|u\|^2,
 \quad \delta > 0 .
\] 
In view of (\ref{4}), we put
\begin{equation}\label{6}
 (u(0),v) = (u^0, v), 
 \quad \forall \ v \in H^1_0(\Omega) .
\end{equation}
The variational (weak) formulation of the problem (\ref{1})--(\ref{4}) consists in finding
$u(\bm x, t) \in H^1_0(\Omega)$, $0 < t \leq T$ that satisfies (\ref{5}), (\ref{6})
with the condition (\ref{2}) on the boundary. 

For the solution of the problem (\ref{5}), (\ref{6}), we have the a priori estimate
\begin{equation}\label{7}
 \|u(t)\| \leq \|u^0\| + \int_{0}^{t}\|f(\theta)\| d \theta .
\end{equation}
To show this, we put $v = u$ in (\ref{5}) and get
\[
 \|u\| \frac{d }{d t}\|u\| + a(u,u) = (f, u ).
\]
Taking into account the positive definiteness of the form $a(\cdot, \cdot)$ and the inequality
\[
 (f,u) \leq \|f\| \|u\| ,
\] 
we obtain
\[
 \frac{d }{d t} \|u(t)\| \leq \|f(t)\| .
\]
From this inequality, in view of (\ref{6}), it follows that the estimate (\ref{7}) holds.

It is easy to obtain less trivial and more interesting a priori estimates for the solution of the problem (\ref{5}), (\ref{6}). We confine ourselves to the elementary estimate (\ref{7}) aiming at obtaining similar estimates using various approximations in time and in space for the problem under consideration.

For numerical solving the IBV problem (\ref{3}), (\ref{4}), approximation in space is constructed using  the finite element method. 
The weak formulation (\ref{5}), (\ref{6}) is employed.
Define the subspace of finite elements $V^h \subset H_0^1(\Omega)$ and the discrete elliptic operator $A$ as
\[
(A y, v) = a(y,v),
\quad \forall \ y,v \in V^h . 
\]
The operator $A$ acts on the finite dimensional space $V^h$ 
and
\begin{equation}\label{8}
A = A^* \geq \delta_h I ,
\quad \delta_h > 0 . 
\end{equation} 
where $I$ is the identity operator.

For the problem (\ref{3}), (\ref{4}), we put into the correspondence the operator equation for $w(t) \in V^h$:
\begin{equation}\label{9}
 \frac{d w}{d t} + A w = \varphi(t), 
 \quad 0 < t \leq T, 
\end{equation} 
\begin{equation}\label{10}
 w(0) = w^0, 
\end{equation} 
where $\varphi(t) = P f(t)$, $w^0 = P u^0$ with $P$ denoting $L_2$-projection onto $V^h$.

Let us multiply equation (\ref{9}) by $w$ and integrate it over the domain $\Omega$:
\[
 \left (\frac{d w}{d t}, w \right ) + (A w, w ) = (\varphi, w) .
\]
In view of the self-adjointness and positive definiteness of the operator $A$, 
we have
\[
 \left (\frac{d w}{d t}, w \right ) \leq  (\varphi, w) .
\]
The right-hand side can be evaluated by the inequality
\[
 (\varphi, w ) \leq  \|\varphi\| \| w \|.
\] 
By virtue of this, we have
\[
 \frac{d}{d t} \|w\| \leq \|\varphi \| .
\] 
The latter inequality leads us to the desired a priori estimate:
\begin{equation}\label{11}
 \|w(t)\| \leq \|w^0\| + \int_{0}^{t}\|\varphi(\theta) \| d \theta .
\end{equation} 
The estimate (\ref{11}) expresses the stability of the solution of the problem (\ref{9}), (\ref{10})
with respect to the initial data and the right-hand side and is quite similar to the estimate (\ref{7}).

To solve numerically the problem (\ref{9}), (\ref{10}),
we apply the implicit two-level scheme \citep{Samarskii1989}.
Let $\tau$ be a step of the uniform grid in time so that
$y^n = y(t^n), \ t^n = n \tau$, $n = 0,1, ..., N, \ N\tau = T$.
For a constant weight parameter $\sigma \ (0 < \sigma \leq 1)$, we define
\[
 y^{n+\sigma} = \sigma y^{n+1} + (1-\sigma) y^{n} .
\]
To approximate equation (\ref{9}), let us consider the following two-level scheme:
\begin{equation}\label{12}
 \frac{y^{n+1} - y^{n}}{\tau} + A y^{n+\sigma} = 
 \varphi^{n+\sigma}, 
  \quad n=0,1,... , N-1,
\end{equation} 
\begin{equation}\label{13}
  y^0 = w^0 .
\end{equation} 
For $\sigma = 0$, the scheme (\ref{12}), (\ref{13}) becomes the explicit scheme, for 
$\sigma = 1$, we get the fully implicit scheme, whereas $\sigma = 0.5$ corresponds to the symmetric scheme (the Crank-Nicolson scheme).
The condition
\[
 I + \tau \left ( \sigma - \frac{1}{2} \right ) A \geq 0 
\]
is necessary and sufficient for the stability of the scheme \citep{Samarskii1989,SamarskiiMatusVabischevich2002}.
In particular, the following statement is true.

\begin{thm}\label{t-1}
The difference scheme (\ref{12}), (\ref{13}) for $\sigma \geq 0.5$ and (\ref{8}) 
is unconditionally stable in $H$ and its solution satisfies the a priori estimate
\begin{equation}\label{14}
 \|y^{n+1}\| \leq \|w^{0}\| + \tau \sum_{j=0}^{n} \|\varphi^{j+\sigma}\|,
  \quad n = 0, 1, ..., N-1 .
\end{equation}
\end{thm} 
 
\begin{pf}
Let us estimate the transition operator. Rewrite the scheme with weights (\ref{12}) in the form
\begin{equation}\label{15}
 y^{n+1} = S y^{n} + \tau  \varphi^{n+\sigma},
\end{equation} 
where $S$ is the operator of transition to a new time level:
\[
 S = S^* = I - \tau B^{-1} A,
 \quad B = I + \sigma \tau A.
\]
Let us formulate the restrictions on the weight $\sigma$ that guarantee the following two-sided inequality:
\begin{equation}\label{16}
 - I \leq S \leq I ,
\end{equation} 
where $\|S\| \leq 1$. 
Taking into account the commutativity of the operators $B$ and $S$,
(\ref{16}) is equivalent to
\[
 - B \leq BS \leq B .
\]
The right inequality is fulfilled for all $\sigma \geq 0$.
The left inequality gives
\[
 2 I + (2 \sigma -1) \tau A \geq 0 .
\]
For non-negative operators $A$, it always holds for $\sigma \geq 0.5$.
For (\ref{15}), in view of (\ref{16}), we get
\[
 \|y^{n+1}\| \leq  \|y^{n}\| + \tau \|\varphi^{n+\sigma}\| ,
\] 
which provides the estimate (\ref{14}).
\end{pf} 

\section{Two-component scheme based on a partition of unity} 

For the differential problem under the consideration, we select the domain decomposition
\begin{equation}\label{17}
  \overline{\Omega} = \overline{\Omega}_{1} + \overline{\Omega}_{2} ,
  \quad \overline{\Omega}_{\alpha} = 
  \Omega_{\alpha} \cup \partial \Omega_{\alpha},
  \quad \alpha = 1, 2
\end{equation}
with overlapping  subdomains
($\Omega_{1 2} \equiv  \Omega_{1 } \cap \Omega_{2} \neq \varnothing$)
\citep{QuarteroniValli1999,ToselliWidlund2005}.
To organize parallel computations, each subdomain consists of a set of disconnected subdomains.

To construct schemes of domain decomposition, we use a partition of unity for the computational domain
$\Omega$  \citep{Laevsky1987,mathew1998domain}. 
Each separate subdomain $\Omega_{\alpha}, \ \alpha = 1,2$
we associate with the function $\eta_{\alpha}({\bm x}), \ \alpha = 1,2$ such that
\[
  \eta_{\alpha}({\bm x}) = \left \{
   \begin{array}{cc}
     > 0, &  {\bm x} \in \Omega_{\alpha},\\
     0, &  {\bm x} \notin  \Omega_{\alpha}, \\
   \end{array}
  \right .
  \quad \alpha = 1,2 ,  
\]
where
\[
  \eta_{1}({\bm x}) + \eta_{2}({\bm x})= 1,
  \quad {\bm x} \in \Omega .
\]
The standard approach \citep{SamarskiiMatusVabischevich2002} is based on  the additive representation
of the operator of the problem (\ref{9}), (\ref{10}):
\begin{equation}\label{18}
  A = A_1 + A_2, 
\end{equation} 
where each individual operator term $A_{\alpha}$ is associated with the separate subdomain 
$\Omega_{\alpha}, \ \alpha = 1,2$.
For instance, in view of (\ref{1}), it is natural to put
\[
 {\cal A}_{\alpha} =  - \nabla (\eta_{\alpha}({\bm x}) k(\bm x) \nabla u ) + \eta_{\alpha}({\bm x}) c(\bm x) u ,
 \quad \alpha = 1,2 ,   
 \quad \bm x \in \Omega . 
\]
In this case, for the representation (\ref{18}), we have
\[
 A_{\alpha}  = A_{\alpha}^* \geq 0,
 \quad \alpha = 1,2 . 
\]  

The construction and investigation of domain decomposition schemes for the unsteady problems
(\ref{9}), (\ref{10}), (\ref{18}) involves the consideration of the appropriate splitting schemes
\citep{Vabishchevich2014}. In the case of two-component splitting, we can apply
the following additive operator--difference schemes of ADI type:
the Douglas--Rachford or Peaceman--Rachford scheme and factorized schemes, which generalize them.

In the factorized scheme, an approximate solution at a new time level is evaluated from the equation
\begin{equation}\label{19}
 B_1 B_2 \frac{y^{n+1} - y^{n}}{\tau} + (A_1 + A_2) y^{n} = \varphi^{n+\sigma} ,
\end{equation} 
where
\begin{equation}\label{20}
 B_\alpha = I + \sigma \tau A_\alpha ,
 \quad \alpha = 1,2 . 
\end{equation} 
For $\sigma = 0.5$, the factorized scheme (\ref{19}), (\ref{20}) 
corresponds to Peaceman--Rachford scheme, whereas  $\sigma = 1$ results in 
the Douglas--Rachford scheme.

For the numerical implementation of the factorized scheme, we can introduce
the auxiliary value $\widetilde{y}^{n+\sigma}$, which is determined from the equation
\begin{equation}\label{21}
 \frac{\widetilde{y}^{n+\sigma} - y^{n}}{\tau} + A_1 \widetilde{y}^{n+\sigma} + A_2 y^{n} = \varphi^{n+\sigma} .
\end{equation}
For the approximate solution at a new  time level, we have
\begin{equation}\label{22}
 \frac{y^{n+1} - y^{n}}{\tau} + A_1 \widetilde{y}^{n+\sigma} + A_2 y^{n+1}  = \varphi^{n+\sigma} .
\end{equation}
Thus, explicit-implicit approximations are used for $A w$ in (\ref{9})
taking into account the splitting (\ref{18}).

Now we provide an estimate for the stability of the factorized scheme, which is similar to the estimate
(\ref{14}) for the scheme with weights (\ref{12}), (\ref{13}). 

\begin{thm}\label{t-2}
The factorized difference scheme (\ref{13}), (\ref{18})--(\ref{20}) is unconditionally stable
for $\sigma \geq 0.5$ and its solution satisfies the a priori estimate
\begin{equation}\label{23}
 \|B_2 y^{n+1}\| \leq \|B_2 w^{0}\| + \tau \sum_{j=0}^{n} \|\varphi^{j+\sigma}\|,
  \quad n = 0, 1, ..., N-1 .
\end{equation}
\end{thm} 
 
\begin{pf}
The consideration is conducted by analogy with the proof of  theorem \ref{t-1}.
From (\ref{19}), we have
\begin{equation}\label{24}
 B_2 y^{n+1} = S B_2  y^{n} + \tau  B_1^{-1} \varphi^{n+\sigma} .
\end{equation} 
In the case of (\ref{19}), (\ref{20}),  the transition operator can be written
\[
 S = \frac{2\sigma -1}{2\sigma}  I + \frac{1}{2\sigma} S_1 S_2 ,
\] 
where
\[
 S_\alpha = (I + \sigma \tau A_\alpha)^{-1}  (I - \sigma \tau A_\alpha),
 \quad \alpha = 1,2 . 
\] 
Using Kellogg’s lemma (see, e.g., \cite{kellogg1964alternating,GrossmannRoosStynes2007,Vabishchevich2014}),
we obtain
\[
 \|S_\alpha\| \leq 1,
 \quad \alpha = 1,2 . 
\]  
For $\sigma \geq 0.5$, we have
\[
 \|S\| \leq \frac{2\sigma -1}{2\sigma}  + \frac{1}{2\sigma}\|S_1\| \|S_2\| \leq 1 .
\] 
In view of this, from (\ref{24}), we get
\[
  \|B_2 y^{n+1}\| \leq  \|B_2  y^{n}\| + \tau \|\varphi^{n+\sigma}\| .
\]
This inequality leads to the estimate (\ref{23}).
\end{pf}

\section{Domain decomposition scheme based on indicator functions of subdomains} 

A new variant of domain decomposition schemes is designed using indicator functions 
for the subdomains $\Omega_\alpha, \ \alpha = 1,2$. Let us introduce
\[
  \chi_{\alpha}({\bm x}) = \left \{
   \begin{array}{cc}
     1, &  {\bm x} \in \Omega_{\alpha},\\
     0, &  {\bm x} \notin  \Omega_{\alpha}, \\
   \end{array}
  \right .
  \quad \alpha = 1,2 . 
\]
Define also the indicator function for the domain of overlap $\Omega_{12}$:
\[
  \chi_{12}({\bm x}) = \left \{
   \begin{array}{cc}
     1, &  {\bm x} \in \Omega_{12},\\
     0, &  {\bm x} \notin  \Omega_{12}. \\
   \end{array}
  \right .
\]
Thus, we have
\begin{equation}\label{25}
 \chi_{1}({\bm x}) + \chi_{2}({\bm x}) - \chi_{12}({\bm x}) = 1 ,
  \quad {\bm x} \in \Omega . 
\end{equation} 

Operators of decomposition are constructed using the indicator functions 
of the subdomains $\Omega_\alpha, \ \alpha = 1,2$ and $\Omega_{12}$:
\[
\begin{split}
 {\cal A}_{\alpha} & =  - \nabla (\chi_{\alpha}({\bm x}) k(\bm x) \nabla u ) + \chi_{\alpha}({\bm x}) c(\bm x) u ,
 \quad \alpha = 1,2 ,  \\
 {\cal A}_{12} & =  - \nabla (\chi_{12}({\bm x}) k(\bm x) \nabla u ) + \chi_{12}({\bm x}) c(\bm x) u ,
 \quad \bm x \in \Omega . 
\end{split}
\]
For  equality ((\ref{25}), we have a three-component representation of the problem operator:
\begin{equation}\label{26}
 A = A_1 + A_2 - A_{12},
\end{equation} 
where
\[
 A_{\alpha}  = A_{\alpha}^* \geq 0,
 \quad \alpha = 1,2 ,
 \quad A_{12}  = A_{12}^* \geq 0 . 
\]
It is necessary to construct splitting schemes for the problem (\ref{9}), (\ref{10}), (\ref{26}),
where a transition to a new time level is performed, as before, by solving problems for the operators
$A_\alpha, \ \alpha =1,2$.

Taking into account the splitting (\ref{26}), similarly to (\ref{21}), (\ref{22}),
an approximate solution is obtained from equations
\begin{equation}\label{27}
 \frac{\widetilde{y}^{n+\sigma} - y^{n}}{\tau} + A_1 \widetilde{y}^{n+\sigma} + A_2 y^{n} 
 -  A_{12} y^{n} = \varphi^{n+\sigma} ,
\end{equation}  
\begin{equation}\label{28}
 \frac{y^{n+1} - y^{n}}{\tau} + A_1 \widetilde{y}^{n+\sigma} + A_2 y^{n+1}  
 - A_{12} \widetilde{y}^{n+\sigma} = \varphi^{n+\sigma} .
\end{equation}
To study the stability of the scheme (\ref{27}), (\ref{28}), we employ theorem \ref{t-2}.

Suppose
\[
 \bar{A}_\alpha = A_\alpha - \frac{1}{2} A_{12},
 \quad \alpha = 1,2 .   
\] 
In view of (\ref{26}), we get
\[
 A = \bar{A}_1 + \bar{A}_2,
 \quad \bar{A}_\alpha = \bar{A}_\alpha^* \geq 0,
 \quad \alpha = 1,2 .   
\]
With the notation introduced above, the scheme (\ref{27}), (\ref{28}) is written as
\begin{equation}\label{29}
 G \frac{\widetilde{y}^{n+\sigma} - y^{n}}{\tau} + \bar{A}_1 \widetilde{y}^{n+\sigma} + \bar{A}_2 y^{n} = \varphi^{n+\sigma} ,
\end{equation}  
\begin{equation}\label{30}
 G \frac{y^{n+1} - y^{n}}{\tau} + \bar{A}_1 \widetilde{y}^{n+\sigma} + \bar{A}_2 y^{n+1} = \varphi^{n+\sigma} ,
\end{equation} 
where
\[
 G = I +  \frac{\tau}{2}  A_{12} ,
 \quad G = G^* \geq I .
\] 

Let
\[
 \widetilde{y}^{n+\sigma}= G^{1/2} v^{n+\sigma},
 \quad  y^{n} = G^{1/2} v^{n},
\] 
\[
 \widetilde{A}_\alpha = G^{-1/2} \bar{A}_\alpha G^{-1/2},
 \quad \alpha =1,2. 
\]
This allows to write  (\ref{29}), (\ref{30}) in the form
\begin{equation}\label{31}
 \frac{\widetilde{v}^{n+\sigma} - v^{n}}{\tau} + \widetilde{A}_1 \widetilde{v}^{n+\sigma} + \widetilde{A}_2 v^{n} = \widetilde{\varphi}^{n+\sigma} ,
\end{equation}  
\begin{equation}\label{32}
\frac{v^{n+1} - v^{n}}{\tau} + \widetilde{A}_1 \widetilde{v}^{n+\sigma} + \widetilde{A}_2 v^{n+1} = \widetilde{\varphi}^{n+\sigma} ,
\end{equation} 
where
\[
 \widetilde{A}_\alpha = \widetilde{A}_\alpha^* \geq 0,
 \quad \alpha = 1,2 ,  
 \quad \widetilde{\varphi}^{n+\sigma} = G^{-1/2} \varphi^{n+\sigma} .
\] 
To study the stability of the scheme (\ref{31}), (\ref{32}), we use theorem \ref{t-2}.
Under the restriction $\sigma \geq 0.5$, the following estimate holds:
\[
 \|\widetilde{B}_2 v^{n+1}\| \leq \|\widetilde{B}_2 v^{0}\| + \tau \sum_{j=0}^{n} \|\widetilde{\varphi}^{j+\sigma}\|,
  \quad n = 0, 1, ..., N-1 ,
\] 
where
\[
 \widetilde{B}_2 = I + \sigma \tau \widetilde{A}_2 .
\] 
This implies
\begin{equation}\label{33}
 \|Q y^{n+1}\| \leq \|Q w^{0}\| + \tau \sum_{j=0}^{n} \|\varphi^{j+\sigma}\|,
  \quad n = 0, 1, ..., N-1 , 
\end{equation} 
where
\[
 Q = \left  ( I + \frac{\tau }{2} A_{12} \right )^{-1} \left ( I + \sigma \tau A_2 - \frac{\sigma -1}{2} \tau A_{12} \right ) .    
\]
The main result of our consideration is the following basic statement on the stability of
schemes of domain decomposition constructed on the basis of the indicator functions of the subdomains.

\begin{thm}\label{t-3}
The splitting scheme (\ref{13}), (\ref{26})--(\ref{28})
is unconditionally stable for $\sigma \geq 0.5$ and its solution satisfies the a priori estimate (\ref{33}).
\end{thm} 

\section{Numerical example} 

Numerical experiments presented here are of comparative nature.
We consider the two-level scheme with weights (\ref{12}), (\ref{13}) as the reference scheme that provides the benchmark
numerical solution for a comparison with the results of  two schemes of domain decomposition  with overlapping subdomains. 
Namely, these are the scheme based on the partition of unity for the domain (\ref{13}), (\ref{18})--(\ref{20})
and the scheme that uses the indicator functions of subdomains (\ref{13}), (\ref{26})--(\ref{28}).
In our consideration, we monitor the proximity of the above decomposition schemes
to the reference scheme. We confine ourselves to the case of schemes with $\sigma =1$. 

We consider the model  problem (\ref{1})--(\ref{4}) in the unit square ($\Omega = [0,1]\times [0,1]$), where
\[
 k(\bm x) = 1,
 \quad  c(\bm x) = 0,
 \quad  f(\bm x,t) = x_1 - x_2,
 \quad  u^0(\bm x,t) = 0 ,
\] 
and $T = 0.1$.

Let $\bar{y}$ be the benchmark numerical solution obtained using the implicit scheme (\ref{12}), (\ref{13}),
whereas $y_1, y_2$ are the solutions from the decomposition scheme (\ref{13}), (\ref{18})--(\ref{20}) and
(\ref{13}), (\ref{26})--(\ref{28}), respectively.
We evaluate the deviation of the solution obtained using the domain decomposition methods
from the benchmark solution derived without domain decomposition as the error:
\[
 \varepsilon_\beta (t) = \|y_\beta - \bar{y}\|,
 \quad \beta = 1,2 . 
\] 

\begin{figure}[tbp]
  \begin{center}
    \includegraphics[width=0.9\linewidth] {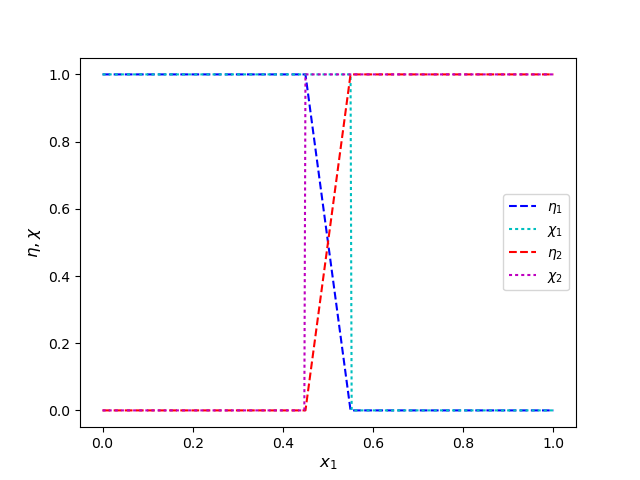}
	\caption{Functions of domain decomposition.}
	\label{fig:1}
  \end{center}
\end{figure}

The computational domain is divided by the variable $x_1$ into two subdomains with the overlap width $2\delta$. 
For the basic variant $\delta=0.05$, the functions $\eta_\alpha(x_1), \ \chi_\alpha(x_1), \ \alpha = 1,2$, 
which define the partition operators, are shown in Fig.\ref{fig:1}. 
The uniform spatial grid $51 \times 51$ was used with piecewise-linear finite elements on triangles.
The number of steps in time is $N = 50$.

\begin{figure}[tbp]
  \begin{center}
    \includegraphics[width=0.9\linewidth] {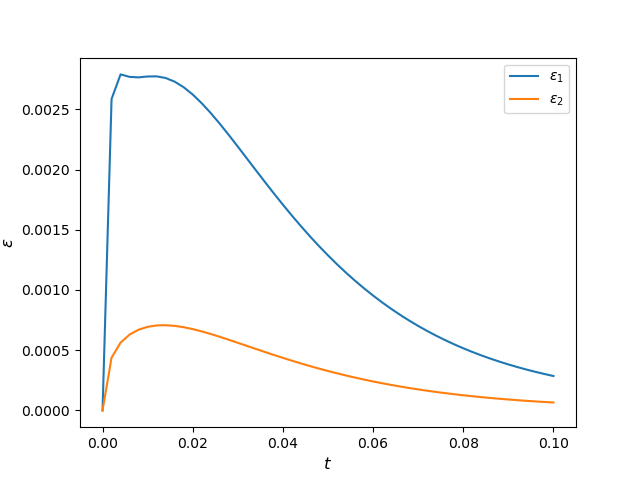}
	\caption{Errors of the domain decomposition schemes.}
	\label{fig:2}
  \end{center}
\end{figure}

\begin{figure}[tbp]
  \begin{center}
    \includegraphics[width=0.9\linewidth] {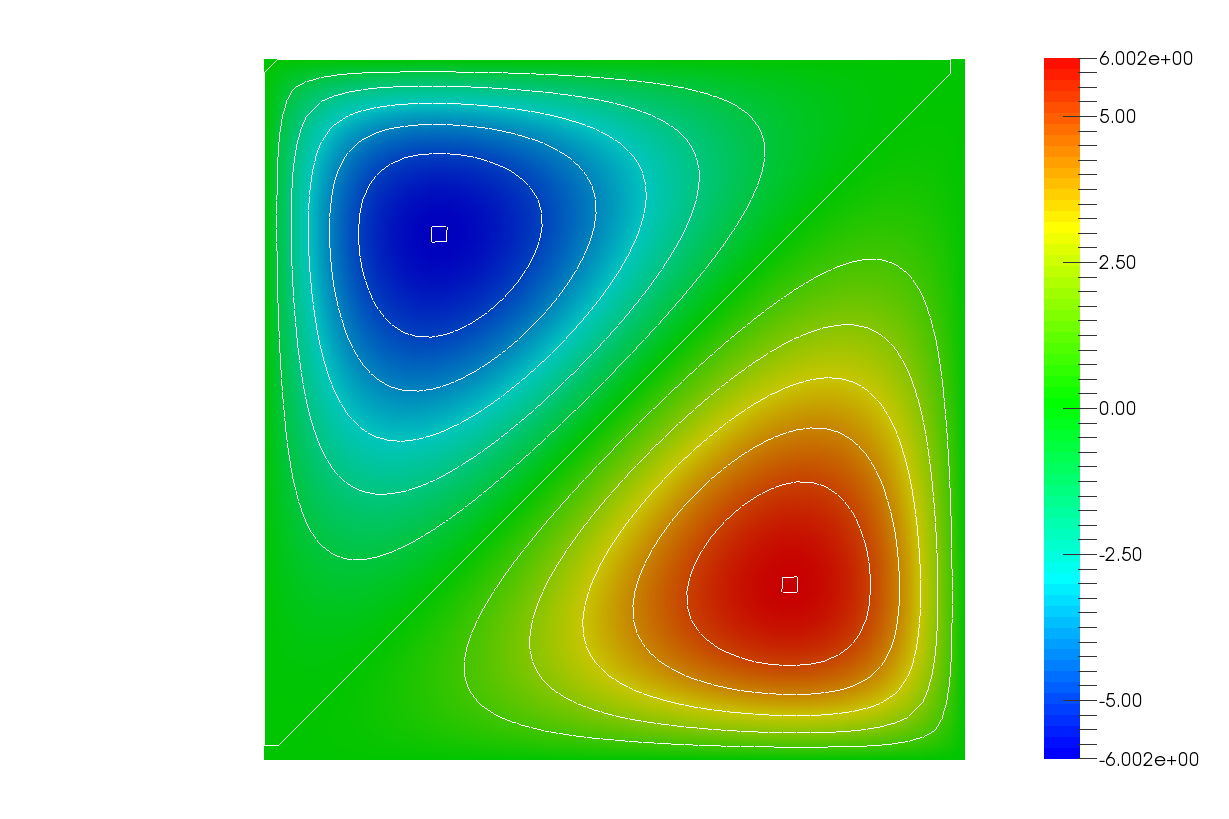}
	\caption{Benchmark solution at the time moment $t=T$.}
	\label{fig:3}
  \end{center}
\end{figure}

\begin{figure}[tbp]
  \begin{center}
    \includegraphics[width=0.9\linewidth] {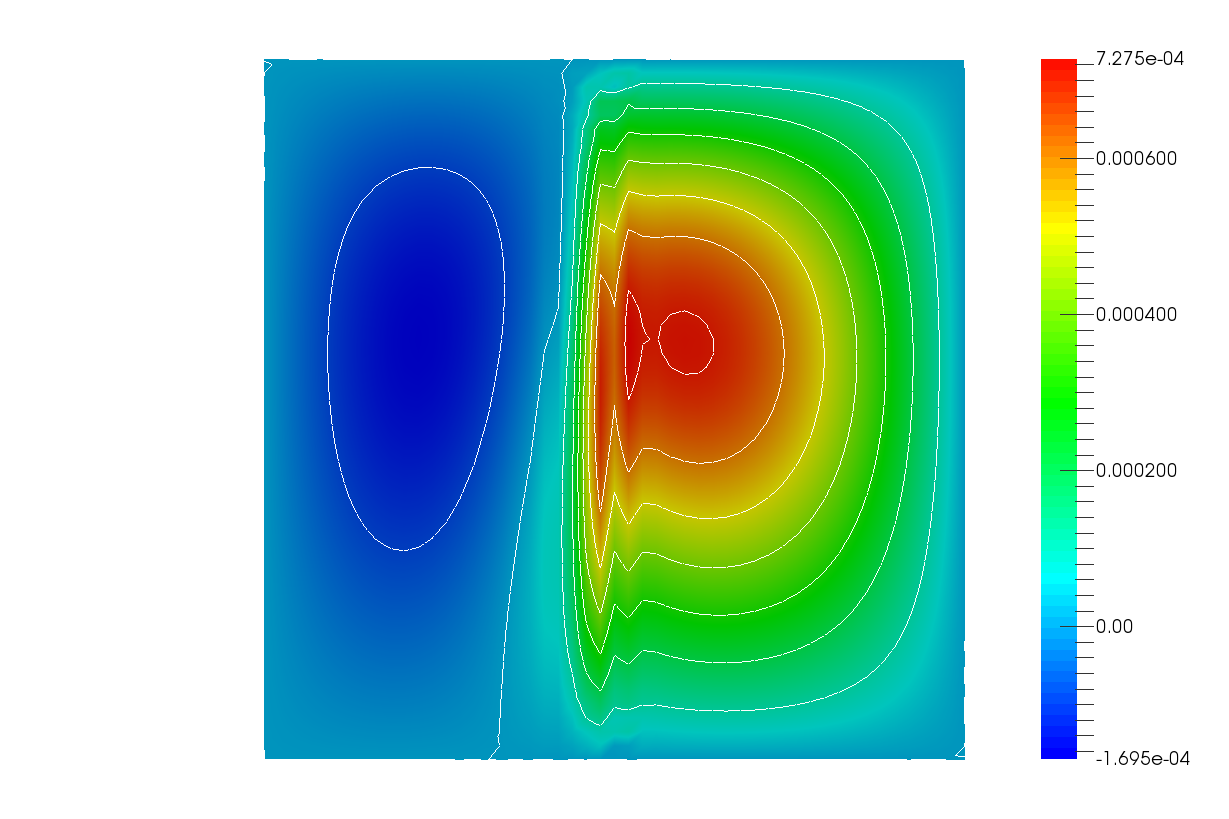}
	\caption{Solution deviation of  the scheme based on the partition of unity.}
	\label{fig:4}
  \end{center}
\end{figure}

\begin{figure}[tbp]
  \begin{center}
    \includegraphics[width=0.9\linewidth] {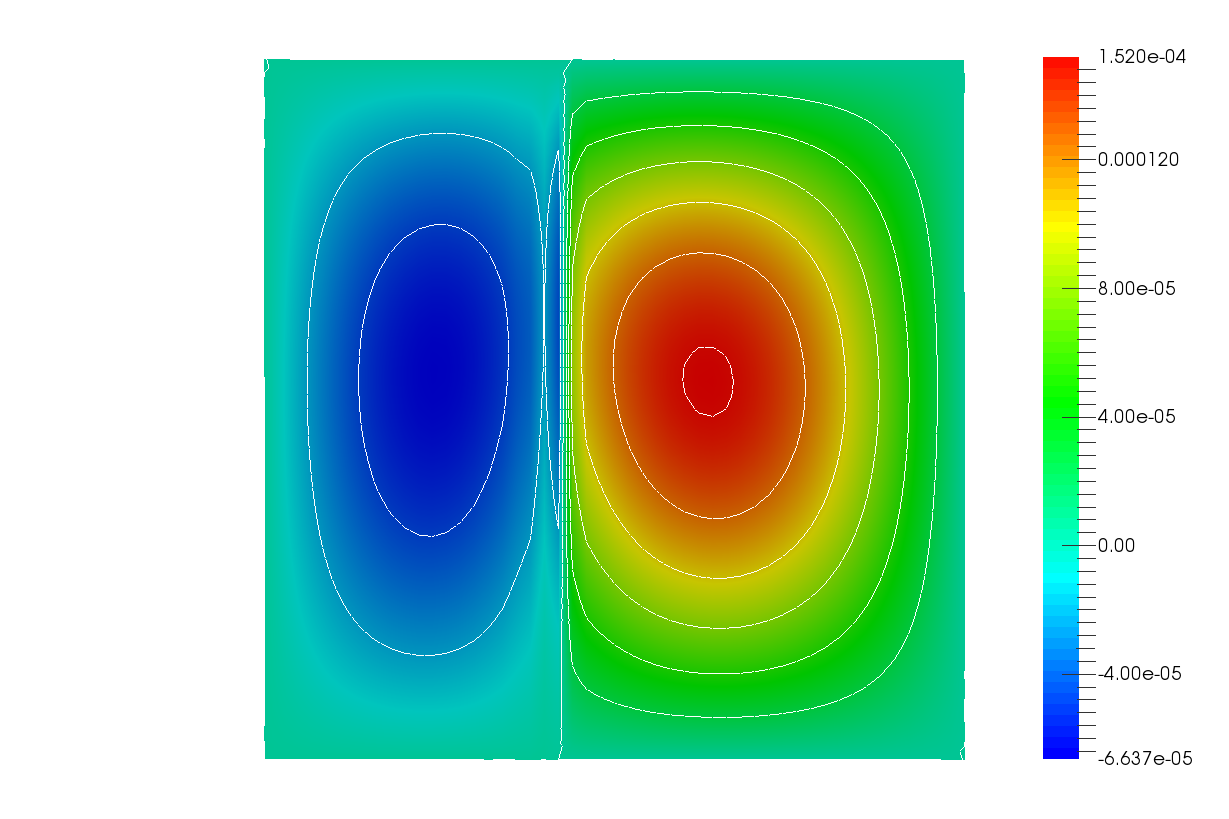}
	\caption{Solution deviation  of the scheme based on the indicator functions.}
	\label{fig:5}
  \end{center}
\end{figure}

The error of the domain decomposition schemes is presented in Fig.\ref{fig:2}. 
In the example considered here, the accuracy of the scheme based on the indicator functions of the subdomains (\ref{13}), (\ref{26})--(\ref{28}) is higher. 
The benchmark solution $\bar{y}$ at the final time moment is shown in  Fig.\ref{fig:3}. The deviation of the solutions obtained with two schemes of domain decomposition 
from the benchmark one is given in Fig.\ref{fig:4},~\ref{fig:5}.

\begin{figure}[tbp]
  \begin{center}
    \includegraphics[width=0.9\linewidth] {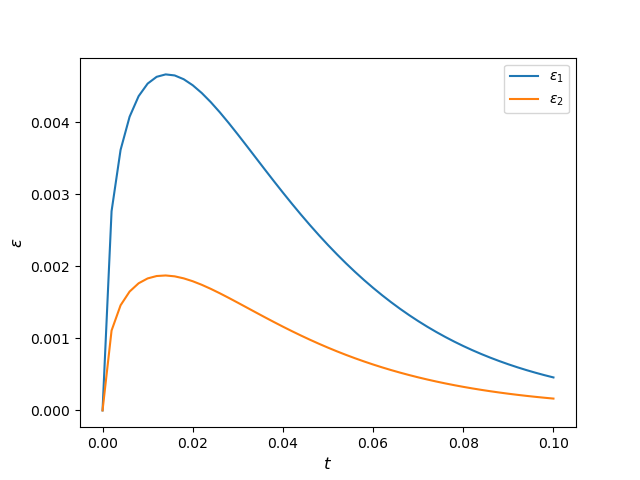}
	\caption{Error of the domain decomposition schemes for  $\delta = 0.025$.}
	\label{fig:6}
  \end{center}
\end{figure}

Reducing the overlap width by half  (see Fig.\ref{fig:6})
results in decreasing the accuracy of the approximate solution.
For $\delta = 0$, both schemes under the consideration coincide and 
lead to the domain decomposition scheme with non-overlapping subdomains.

\begin{figure}[tbp]
  \begin{center}
    \includegraphics[width=0.9\linewidth] {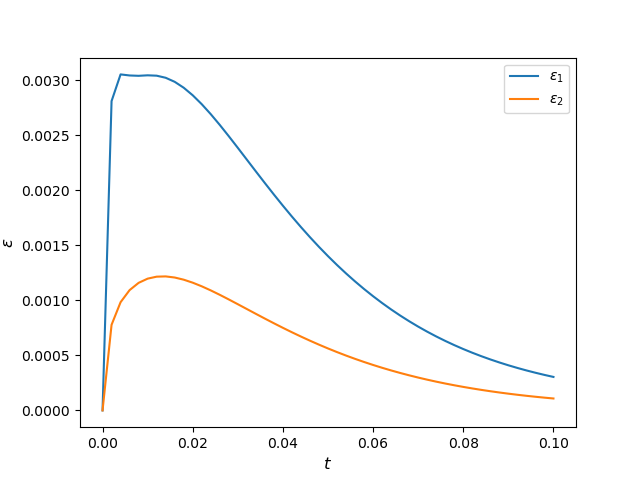}
	\caption{Error of the domain decomposition schemes for the grid  $101 \times 101$.}
	\label{fig:7}
  \end{center}
\end{figure}

\begin{figure}[tbp]
  \begin{center}
    \includegraphics[width=0.9\linewidth] {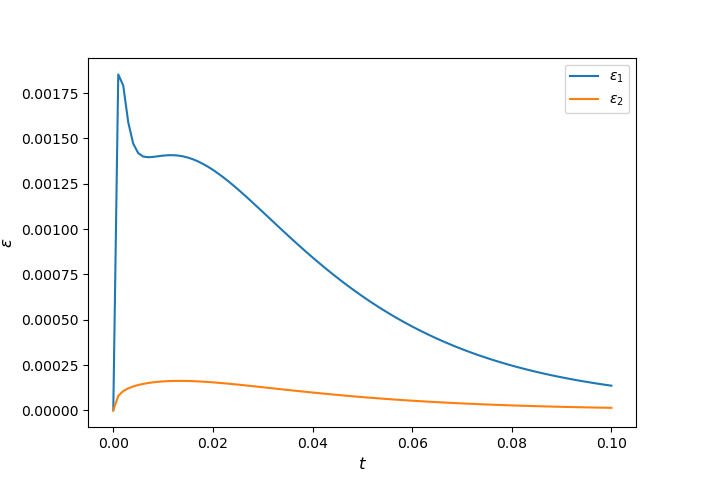}
	\caption{Error of the domain decomposition schemes for $N = 100$.}
	\label{fig:8}
  \end{center}
\end{figure}

The use of a finer grid in space (compare Fig.\ref{fig:2} and Fig.\ref{fig:7}) leads to a drop in accuracy.
This fact demonstrates a conditional convergence of the domain decomposition schemes.
Calculations with a smaller time step (see Fig.\ref{fig:8}) allow to increase the accuracy of the numerical solution.

\clearpage

\section{Conclusions} 

The boundary value problem for a second-order parabolic equation is considered.
The standard finite element approximation in space is employed.
For approximation in time, the two-level scheme with weight ($\theta$-method) is used.

Iteration-free schemes of domain decomposition with overlapping subdomains
are constructed on the basis of a partition of unity for a computational domain.
In the case of a two-component domain decomposition, approximation in time is designed using
factorized schemes, which generalize the standard ADI schemes.

A new class of domain decomposition schemes with overlapping subdomains based on indicator functions of subdomains is proposed for unsteady problems.
In this case, for the problem operator, we have a three-component representation with decomposition operators
connected with subdomains and their intersection.
The unconditional stability of these schemes is established for two-component schemes of domain decomposition.

Numerical results for the model two-dimensional problem demonstrate the robustness of the new decomposition scheme with overlapping subdomains. 
A comparison of the conventional scheme based on a partition of unity with the scheme that used the indicator functions of the subdomains demonstrates its advantage in accuracy.

\section*{Acknowledgements}

This work was supported by the grant of Russian Federation Government (agreement \#~14.Y26.31.0013).


\begin{thebibliography}{23}
\expandafter\ifx\csname natexlab\endcsname\relax\def\natexlab#1{#1}\fi
\expandafter\ifx\csname url\endcsname\relax
  \def\url#1{\texttt{#1}}\fi
\expandafter\ifx\csname urlprefix\endcsname\relax\def\urlprefix{URL }\fi

\bibitem[{Abrashin(1990)}]{Abrashin1990}
Abrashin, V.~N., 1990. {A variant of the method of variable directions for the
  solution of multidimensional problems of mathematical-physics. I.}
  Differential equations 26, 314--323, in Russian.

\bibitem[{Abrashin and Vabishchevich(1998)}]{AbrashinVabischevich1998}
Abrashin, V.~N., Vabishchevich, P.~N., 1998. Vector additive schemes for
  evolution equations of second order. Differential equations 34~(12),
  1666--1674, in Russian.

\bibitem[{Dryja(1991)}]{Dryja1991}
Dryja, M., 1991. Substructuring methods for parabolic problems. In: Glowinski,
  R., Kuznetsov, Y.~A., Meurant, G.~A., P{\'e}riaux, J., Widlund, O. (Eds.),
  Fourth International Symposium on Domain Decomposition Methods for Partial
  Differential Equations. SIAM, Philadelphia, PA.

\bibitem[{Grossmann et~al.(2007)Grossmann, Roos, and
  Stynes}]{GrossmannRoosStynes2007}
Grossmann, C., Roos, H.~G., Stynes, M., 2007. Numerical Treatment of Partial
  Differential Equations. Springer Verlag.

\bibitem[{Kellogg(1964)}]{kellogg1964alternating}
Kellogg, R.~B., 1964. An alternating direction method for operator equations.
  Journal of the Society for Industrial and Applied Mathematics 12~(4),
  848--854.

\bibitem[{Kuznetsov(1988)}]{Kuznetsov1988}
Kuznetsov, Y.~A., 1988. New algorithms for approximate realization of implicit
  difference schemes. Sov. J. Numer. Anal. Math. Model. 3~(2), 99--114.

\bibitem[{Kuznetsov(1991)}]{Kuznetsov1991}
Kuznetsov, Y.~A., 1991. Overlapping domain decomposition methods for
  {FE}-problems with elliptic singular perturbed operators. Fourth
  international symposium on domain decomposition methods for partial
  differential equations, Proc. Symp., Moscow/Russ. 1990, 223-241 (1991).

\bibitem[{Laevsky(1987)}]{Laevsky1987}
Laevsky, Y.~M., 1987. Domain decomposition methods for the solution of
  two-dimensional parabolic equations. In: Variational-difference methods in
  problems of numerical analysis. No.~2. Comp. Cent. Sib. Branch, USSR Acad.
  Sci., Novosibirsk, pp. 112--128, in Russian.

\bibitem[{Marchuk(1990)}]{Marchuk1990}
Marchuk, G.~I., 1990. Splitting and alternating direction methods. In: Ciarlet,
  P.~G., Lions, J.-L. (Eds.), Handbook of Numerical Analysis, Vol. I.
  North-Holland, pp. 197--462.

\bibitem[{Mathew et~al.(1998)Mathew, Polyakov, Russo, and
  Wang}]{mathew1998domain}
Mathew, T.~P., Polyakov, P.~L., Russo, G., Wang, J., 1998. Domain decomposition
  operator splittings for the solution of parabolic equations. SIAM Journal on
  Scientific Computing 19~(3), 912--932.

\bibitem[{Quarteroni and Valli(1999)}]{QuarteroniValli1999}
Quarteroni, A., Valli, A., 1999. Domain Decomposition Methods for Partial
  Differential Equations. Clarendon Press.

\bibitem[{Samarskii(1967)}]{Samarskii1967}
Samarskii, A.~A., 1967. Regularization of difference schemes. Zh. Vychisl. Mat.
  Mat. Fiz. 7, 62--93, in Russian.

\bibitem[{Samarskii(2001)}]{Samarskii1989}
Samarskii, A.~A., 2001. The Theory of Difference Schemes. Marcel Dekker, New
  York.

\bibitem[{Samarskii et~al.(2002)Samarskii, Matus, and
  Vabishchevich}]{SamarskiiMatusVabischevich2002}
Samarskii, A.~A., Matus, P.~P., Vabishchevich, P.~N., 2002. Difference Schemes
  with Operator Factors. Kluwer Academic Pub.

\bibitem[{Samarskii and Vabishchevich(1995)}]{SamarskiiVabishchevich1995a}
Samarskii, A.~A., Vabishchevich, P.~N., 1995. Vector additive schemes of domain
  decomposition for parabolic problems. Differential equations 31, 1563--1569,
  in Russian.

\bibitem[{Toselli and Widlund(2005)}]{ToselliWidlund2005}
Toselli, A., Widlund, O., 2005. Domain Decomposition Methods -- Algorithms and
  Theory. Springer.

\bibitem[{Vabishchevich(1989)}]{Vabischevich1989}
Vabishchevich, P.~N., 1989. Difference schemes decompose the computational
  domain for solving transient problems. Zh. Vychisl. Mat. Mat. Fiz. 29~(12),
  1822--1829, in Russian.

\bibitem[{Vabishchevich(1994{\natexlab{a}})}]{Vabischevich1994b}
Vabishchevich, P.~N., 1994{\natexlab{a}}. Parallel domain decomposition
  algorithms for time-dependent problems of mathematical physics. In: Advances
  in Numerical Methods and Applications. World Schientific, pp. 293--299.

\bibitem[{Vabishchevich(1994{\natexlab{b}})}]{Vabischevich1994a}
Vabishchevich, P.~N., 1994{\natexlab{b}}. Regionally additive difference
  schemes for stabilizing correction for parabolic problems. Zh. Vychisl. Mat.
  Mat. Fiz. 34~(12), 1832--1842, in Russian.

\bibitem[{Vabishchevich(1996)}]{Vabischevich1996b}
Vabishchevich, P.~N., 1996. Vector additive difference schemes for second order
  evolution equations. Zh. Vychisl. Mat. Mat. Fiz. 36(3), 44--51, in Russian.

\bibitem[{Vabishchevich(2008)}]{Vabishchevich2008}
Vabishchevich, P.~N., 2008. Domain decomposition methods with overlapping
  subdomains for the time-dependent problems of mathematical physics.
  Computational Methods in Applied Mathematics 8~(4), 393--405.

\bibitem[{Vabishchevich(2011)}]{vab-2011}
Vabishchevich, P.~N., 2011. A substructuring domain decomposition scheme for
  unsteady problems. Computational Methods in Applied Mathematics 11~(2),
  241--268.

\bibitem[{Vabishchevich(2014)}]{Vabishchevich2014}
Vabishchevich, P.~N., 2014. {Additive Operator-Difference Schemes: Splitting
  Schemes}. de Gruyter.

\end{thebibliography}
\end{document}